\documentclass[leqno,twoside, 11pt]{amsart}
\usepackage{amssymb, latexsym, esint}
\usepackage{color}
\usepackage{amsmath,amsthm,amscd,amsfonts,mathrsfs}
\usepackage{enumerate}

\theoremstyle{plain}
\numberwithin{equation}{section} \numberwithin{figure}{section}

\newtheorem{theorem}{Theorem}[section]
\newtheorem{lemma}[theorem]{Lemma}
\newtheorem{proposition}[theorem]{Proposition}
\newtheorem{corollary}[theorem]{Corollary}
\newtheorem{defn}[theorem]{Definition}
  \usepackage{epsfig} 
\usepackage{graphicx}  \usepackage{epstopdf}
\theoremstyle{definition}
\newtheorem{remark}[theorem]{Remark}

\begin{document}
\markboth{Pablo Ochoa and Virginia N. Vera de Serio }
{}

%

%
\title[Viscosity solutions and optimal value function]{A viscosity solution approach to regularity properties of the optimal value function}

%
%

\author{Pablo Ochoa}

\address{Pablo Ochoa. Universidad Nacional de Cuyo, Facultad de Ingenieri\'ia-CONICET\\Parque Gral. San Mart\'in 5500\\
Mendoza, Argentina.}
\email{pablo.ochoa@ingenieria.uncuyo.edu.ar}

\author{Virginia N. Vera de Serio}%

\address{Virginia N. Vera de Serio. Universidad Nacional de Cuyo, Facultad de Ciencias Econ\'omicas, Mendoza, Argentina.}%
\email{virginia.vera@fce.uncu.edu.ar}


%



\maketitle

\begin{abstract}%
In this paper we analyze the optimal value function $v$ associated to a
general parametric optimization problems via the theory of viscosity
solutions. The novelty is that we obtain regularity properties of $v$ by
showing that it is a viscosity solution to a set of first-order equations.
As a consequence, in Banach spaces, we provide sufficient conditions for
local and global Lipschitz properties of $v$. We also derive, in finite
dimensions, conditions for optimality through a comparison principle.
Finally, we study the relationship between viscosity and Clarke generalized
solutions to get further differentiability properties of $v$ in Euclidean
spaces.

\end{abstract}%

%



\section{Introduction}

We consider parametric optimization problems of the form 
\begin{equation}
\left( P_{u}\right) :\inf_{x\in \Phi (u)}f(x,u)  \label{problem}
\end{equation}%
depending on a parameter vector $u\in U.$ Here $x$ belongs to a non empty
subset $X$ of a real Banach space $\mathfrak{X,}$ $U$\ is a non empty subset
of a real Banach space $\mathcal{Y}$, $f:X\mathfrak{\times }U\mathfrak{%
\rightarrow \mathbb{R}}$ is a continuous function, and $\Phi :U\rightarrow
2^{X}$ is the feasible set mapping. The corresponding optimal value function 
$v:U\rightarrow \left[ -\infty ,+\infty \right] $ is defined as 
\begin{equation*}
v(u):=\inf_{x\in \Phi (u)}f(x,u),\quad u\in U,
\end{equation*}%
and the optimal set function $S:U\rightarrow 2^{X}$ is given by 
\begin{equation*}
S\left( u\right) :=\left\{ x\in \Phi \left( u\right) :v\left( u\right)
=f\left( x,u\right) \right\} .
\end{equation*}

Our main regularity result concerns Lipschitz and differentiability
properties\ of $v$ that are obtained by appealing to the theory of viscosity
solutions. We study the Lipschitz continuity of $v$ for Banach spaces in two
steps. First, we provide conditions related to certain first-order partial
differential equation that imply this Lipschitz property, via the viscosity
solution concept. Second, we analyze conditions based on the data that give
the validity of these sufficient conditions for $v$.

The Lipschitzian properties of the optimal value function and of the optimal
set function have been studied by a large number of authors since the 80's.
Indeed, Aubin in 1984, \cite{Aub}, addressed the Lipschitz property of the
optimal set function $S,$ in Banach spaces. Since then, we can cite, among
others, \cite{R85}, \cite{GJ89}, \cite{M97}, \cite{CKLP}, and \cite{DKV} in $%
\mathbb{R}^{n},$ and \cite{CDLP}, \cite{DM}, and \cite{DR} for Banach
spaces, that have also analyzed this regularity property of $S$ in different
situations, e.g. linear, convex, non-linear cases, optimal control problems,
semi-infinite programming. With respect to the Lipschitz properties of $v,$
we can refer to \cite{B84}, \cite{CLPT99}, \cite{CLPT2006}, \cite{T2008},
and \cite{GT2019} in $\mathbb{R}^{n};$ \cite{DM}, and \cite{GH2017} for
Banach spaces. See the surveys \cite{L2012}\ and \cite{GL2014}, and
references therein, for Lipschitz properties of $v$ and $S$ in the case of
semi-infinite optimization. The directional lipschitzian optimal solutions
and directional derivative for the optimal value function have also been
studied, see, e.g., \cite{R85}, \cite{GJ89}, and \cite{T91}. We refer to 
\cite{BS} and references therein for a throughout presentation of
perturbation analysis of optimization problems for Banach spaces.

In this paper we propose a new approach via the viscosity solution concept.
Viscosity solutions were introduced by Crandall and Lions in \cite{CL} for
first-order differential equations of Hamilton-Jacobi type. In general,
nonlinear first-order problems do not have classical solutions, even if the
Hamiltonian and boundary conditions are smooth. Indeed, these kind of
equations are usually satisfied by value functions, which fail to be $C^{1}$
or even differentiable. Value functions arise in deterministic control
problems, deterministic differential games and the calculus of variations.
The book \cite{BCD} contains a nice presentation of the relationship between
value functions and viscosity solutions, and further aspects of the theory
of Hamilton-Jacobi-Bellman equations. We point out that the advantage of the
weak type of solution introduced in \cite{CL} over other notions is that it
provides existence (\cite{CL}, \cite{I}, \cite{Sou}), stability (\cite{CEL}, 
\cite{Sou}) and uniqueness (\cite{CL}, \cite{CEL}) of solutions. Regularity
issues, like Lipschitz and H\"{o}lder continuity, $C^{1}$ regularity and
second-order differentiability may also be obtained for viscosity solutions
of first-order equations (\cite{BCD}, \cite{I}). The theory was also
extended to second-order degenerate elliptic equations in \cite{L83} and 
\cite{L83b}, and to infinite dimensional spaces in \cite{CLvarios}. We
highlight the reference \cite{J} where a proof of uniqueness by purely PDE
methods is provided. A review of the theory of viscosity solutions, with a
vast source of references, can be found in \cite{CIL}.

Our assay resembles the analysis in \cite{CF2017} referred to the regularity
properties of the value function of infinite horizon problem in optimal
control, except that instead of using generalized gradients, we mostly
consider super and subgradients in viscosity sense. In particular, the paper 
\cite{F} that discusses the relationship of locally Lipschitz continuous
viscosity solutions and generalized derivative in Clarke sense for
Hamilton-Jacobi equations is very close to the final section of the present
work. One main difference with the application to optimal control theory is
that we postulate a new first-order equation to tackle the problem (\ref%
{problem}) and find conditions on the data so that a viscosity solution is
locally Lipschitz.

Furthermore, in the context of Euclidean spaces, we obtain a Comparison
Principle that gives the uniqueness property and allows a verification
technique to determine sufficient conditions for both, optimality and non
optimality. In addition, we provide a necessary condition for optimality. In 
\cite{Z1991}, for optimal control theory, it is developed a verification
technique in the framework of viscosity solutions, here we provide somewhat
analogous conditions for general parametric optimization problems. Finally,
in this finite dimensional setting, we derive relationships between the
viscosity solutions and the Clarke generalized solutions. These relations
give the possibility of obtaining pointwise differentiability properties of $%
v$ and deduce, for convex problems, that $v$ is $C%
{{}^1}%
$ under some conditions on the data $f$ and $\Phi $. This analysis of going
from pointwise differentiability to local differentiability is developed for
Hamilton-Jacobi equations in \cite{CF2013}.

This paper is organized as follows. Section 2 introduces the necessary
notation and preliminaries. In Section 3 we discuss the Lipschitz continuity
of viscosity solutions of a certain first-order equation. Section 4 provides
conditions on the data so that the optimal value function is a viscosity
solution of this equation, and hence Lipschitz continuous. Finally, in
Section 5 we restrict our analysis to finite dimensional spaces to find a
comparison principle, a uniqueness property, and optimality conditions.
Furthermore, we study the relationship between locally Lipschitz continuous
viscosity solutions and solutions in the sense of Clarke generalized
derivative, obtaining differentiability properties of the optimal value
function.

\section{Notation and preliminaries}

In this section, we introduce the basic notation, definitions, and results
that we employ throughout the paper.

\subsection{Basic notation and definitions}

For given real Banach space $\mathfrak{B}$, its dual is $\mathfrak{B}^{\ast
} $ and $\left\langle \cdot ,\cdot \right\rangle $ is the usual pairing. The
underlying norm in $\mathfrak{B}$ will be denoted by $||\cdot ||$ (or $%
||\cdot ||_{\mathfrak{B}}$ when is needed for clarity) and the closed ball
centered at $z\in \mathfrak{B}$ with radius $r>0$ will be $B(z,r)$. Whenever
we consider a Cartesian product of topological spaces we will always endow
it with the product topology.

For $E\subset \mathfrak{B}$, $USC(E),LSC(E)$\ and $C(E)$ will stand for the
sets of\textnormal{\ upper semicontinuous, lower semicontinuous, and
continuous functions }$w:E\rightarrow \mathbb{R}$, respectively. We say that 
$g:\mathfrak{B}\rightarrow \mathbb{R}$ is Fr\'{e}chet differentiable at $%
z\in \mathfrak{B}$ if there is $L_{z}\in \mathfrak{B}^{\ast }$ so that 
\begin{equation*}
\lim_{h\rightarrow 0}\dfrac{g(z+h)-g(z)-\left\langle L_{z},h\right\rangle }{%
\Vert h\Vert }=0,
\end{equation*}%
and we write 
\begin{equation*}
\nabla g(z):=L_{z}\text{ \ \ \ and \ \ \ }\nabla g(z)\cdot h=\left\langle
\nabla g(z),h\right\rangle =\left\langle L_{z},h\right\rangle .
\end{equation*}%
We will mostly use the notation $\nabla g(z)\cdot h$ when the space is of
finite dimension. For open $E\subset \mathfrak{B},C^{1}(E):=\left\{
w:E\rightarrow \mathbb{R},\nabla w\mathnormal{\ }\text{\textnormal{is
continuous in }}E\right\} .$ A norm $\Vert \cdot \Vert $ in $\mathfrak{B}$
is called Fr\'{e}chet differentiable if $\Vert \cdot \Vert $ is Fr\'{e}chet
differentiable at every point of $\mathfrak{B}\setminus \left\{ 0\right\} $.
It is known, {\cite[Theorem 8.19]{FHM},} that\ a Banach space $\mathfrak{B}$
admits an equivalent Fr\'{e}chet differentiable norm whenever $\mathfrak{B}%
^{\ast }$ is separable.

Next we recall some well-known concepts for set-valued mappings. Given two
Banach spaces $\mathfrak{X}$ and $\mathcal{Y}$, and a set $U$ in $\mathcal{Y}
$, a set-valued mapping $\mathcal{A}:U\rightarrow 2^{\mathfrak{X}}$ is said
to be closed at $u_{0}\in U$ if for any sequences $\left\{ x_{k}\right\}
\subset \mathfrak{X},$ $\left\{ u_{k}\right\} \subset U$, such that $%
x_{k}\in \mathcal{A}\left( u_{k}\right) ,x_{k}\rightarrow x_{0},$ and $%
u_{k}\rightarrow u_{0}$, it holds that $x_{0}\in \mathcal{A}\left(
u_{0}\right) $. Moreover, $\mathcal{A}$ is called closed if its graph%
\begin{equation*}
\text{gph}\,\mathcal{A}:=\left\{ (u,x):u\in U,x\in \mathcal{A}(u)\right\}
\end{equation*}%
is closed in $U\times \mathfrak{X}$. Also, $\mathcal{A}$ is said to be lower
semicontinuous at $u_{0}$ (lsc, in brief), in the sense of Berge, if for
each open set $W\subset \mathfrak{X}$ such that $W\cap \mathcal{A}%
(u_{0})\neq \emptyset $, there exists an open neighborhood $V$ of $u_{0}$ in 
$U$, such that $\mathcal{A}(u)\cap W\neq \emptyset $ for all $u\in V$. $%
\mathcal{A}$ is said to be upper semicontinuous at $u_{0}$ (usc, in brief),
in the sense of Berge, if for each open set $W\subset \mathfrak{X}$ such
that $\mathcal{A}(u_{0})\subset W$, there is an open neighborhood $V$ of $%
u_{0}$ in $U$ so that $\mathcal{A}(u)\subset W$ for all $u\in V$.

\begin{defn}[Inf-compactness property]
Given $\mathcal{A}:U\rightarrow 2^{\mathfrak{X}},$ a function $f:\mathfrak{X}%
\times U\rightarrow \mathbb{R}$, satisfies the inf-compactness property on
gph $\mathcal{A}$ if for every $u\in U$, there exist $\alpha \in \mathbb{R}$
and a compact set $C\subset \mathfrak{X}$ so that%
\begin{equation*}
\emptyset \neq \left\{ x\in \mathcal{A}(u^{\prime }):f(x,u^{\prime })\leq
\alpha \right\} \subset C
\end{equation*}%
for all $u^{\prime }$ in a neighborhood of $u$.
\end{defn}

For these functions $f$ and for a fixed $x\in \mathfrak{X,}$ we will write $%
f_{x}\left( u\right) :=f\left( x,u\right) ,$ and use indistinctly the
notation $D_{d}f_{x}\left( u\right) $ or $D_{d}f\left( x,u\right) $ for the
directional derivative of the function $f_{x}\left( \cdot \right) $ at a
point $u\in U$ in the direction of $d\in \mathcal{Y}$. We may also use $%
D_{d}f$, for short.

Finally, the following well-known result accounts for continuity of the
optimal value function (see, e.g., \cite[Theorem 1.2]{FKV}). We state it in
the current framework of Banach spaces, but the result holds in any
compactly generated topological space.

\begin{theorem}[Continuity of $v$ and usc of $S$]
\label{Th. v-continuous}Let the problem (\ref{problem}) and suppose that $%
\Phi $ is lower semicontinuous in $U$ and that $f$ is continuous in $%
\mathfrak{X}\times U$ and satisfies the inf-compactness property on gph $%
\Phi .$ Then, the optimal value function $v$ is continuous in $U$ and the
optimal set function $S$ is upper semicontinuous in $U$.
\end{theorem}

\subsection{Viscosity solutions in Banach spaces}

We now recall some useful concepts and results about viscosity solutions. In
what follows, $F:U\times \mathbb{R}\times \mathcal{Y}^{\ast }\rightarrow 
\mathbb{R}$ is a first-order operator.

\begin{defn}
Let $U$ be an open subset of a Banach space $\mathcal{Y}$. We say that a
function $w\in USC(U)$ is a viscosity subsolution\ of 
\begin{equation}
F(u,w,\nabla w)=0  \label{eqPDE}
\end{equation}%
at $u_{0}\in U$ if for any $\eta $ in the set 
\begin{equation*}
\mathcal{J}^{+}w(u_{0}):=\left\{ p\in \mathcal{Y}^{\ast
}:\limsup_{u\rightarrow u_{0}}\frac{w(u)-w(u_{0})-\left\langle
p,u-u_{0}\right\rangle }{\Vert u-u_{0}\Vert }\leq 0\right\} ,
\end{equation*}%
the following holds 
\begin{equation*}
F(u_{0},w(u_{0}),\eta )\leq 0.
\end{equation*}%
Similarly, we say that a function $w\in LSC(U)$ is a viscosity supersolution
of \eqref{eqPDE} at $u_{0}\in U$ if for any $\eta $ in the set 
\begin{equation*}
\mathcal{J}^{-}w(u_{0}):=\left\{ p\in \mathcal{Y}^{\ast
}:\liminf_{u\rightarrow u_{0}}\frac{w(u)-w(u_{0})-\left\langle
p,u-u_{0}\right\rangle }{\Vert u-u_{0}\Vert }\geq 0\right\} ,
\end{equation*}%
the following holds 
\begin{equation*}
F(u_{0},w(u_{0}),\eta )\geq 0.
\end{equation*}%
A viscosity solution of \eqref{eqPDE} at $u_{0}\in U$ is a viscosity sub-
and supersolution at $u_{0}\in U$. Finally, a function is a viscosity
subsolution (supersolution, solution) of \eqref{eqPDE} in $U$ if it is so
for all $u_{0}\in U.$
\end{defn}

We have the following known equivalence, adapted from Proposition 1 in \cite[%
Part 1]{CLvarios}:

\begin{proposition}
Let $w\in USC(U)$. Then, $w$ is a viscosity subsolution of (\ref{eqPDE}) at $%
u_{0}\in $ $U$ if and only if for every $g\in C(U)$ Fr\'{e}chet
differentiable at $u_{0}$ and such that $w-g$ attains a local maximum at $%
u_{0}$, the following holds\ 
\begin{equation*}
F(u_{0},w(u_{0}),\nabla g(u_{0}))\leq 0
\end{equation*}%
An analogous statement holds true for a supersolution.
\end{proposition}

The main obstacle in the application of this last proposition is that closed
bounded sets are not compact with respect to the strong topology (unless the
underlying space is finite dimensional). Hence continuous functions over
these sets may not attain their maximum. However, when the domain under
consideration satisfies the Radon-Nikodym property, any continuous function
attains its maximum under a small linear perturbation. Recall that a closed,
convex and bounded set $D$ of a Banach space $\mathfrak{X}$ is said to be a
RNP set (a set with the Radon-Nikodym property) if for every finite measure
space $(\Omega ,\Sigma ,\mu )$, every vector measure $m:\Sigma \rightarrow 
\mathfrak{X}$ that is of bounded variation, absolutely continuous with
respect to $\mu $, and has average range 
\begin{equation*}
\left\{ \mu (E)^{-1}m(E):E\in \Sigma ,\mu (E)>0\right\}
\end{equation*}%
contained in $D$ is representable by a Bochner integrable function (see,
e.g., \cite{Di} for further details on RNP properties).

According to \cite{St}, examples of RNP sets are convex weakly compact sets
in Banach spaces. Hence, in a reflexive Banach space, closed balls are RNP
sets. The next known theorem is the key in producing maximum of bounded
semicontinuous functions in bounded sets (see, e.g., \cite[pag. 174-176]{St}
for the proof).

\begin{theorem}
\label{RNP set max} Let $\mathcal{Y}$ be a Banach space and let $D\subset 
\mathcal{Y}$ be a RNP set. Assume that $g:D\rightarrow \mathbb{R}$ is upper
semicontinuous and bounded above in $D$. Then, for any $\delta >0$, there
exists $h\in \mathcal{Y}^{\ast }$, $\Vert h\Vert \leq \delta $, such that $%
g+h$ attains its maximum at a point $x_{0}\in D$.
\end{theorem}

\section{Lipschitz continuity of a viscosity subsolution}

In this section, we consider the following first order equation defined in $%
U $, 
\begin{equation}
-\left\langle \nabla w\left( u\right) ,d\right\rangle +\inf_{x\in S\left(
u\right) }D_{d}f_{x}\left( u\right) =0\text{,}  \label{eq}
\end{equation}%
where we assume that the second term has sense. We will provide conditions
that assure that a viscosity subsolution is locally Lipschitz in $U.$

In order to state our regularity result, we will need the following
assumptions.

\textbf{Assumptions on the Banach space }$\mathcal{Y}$\textbf{:}

(H1) $\mathcal{Y}$ is reflexive;

(H2) There exists a function $M:\mathcal{Y}\times \mathcal{Y}\rightarrow
\lbrack 0,\infty ]$ so that: $y_{1}\rightarrow M(y_{1},y_{2})$ and $%
y_{2}\rightarrow M(y_{1},y_{2})$ are Fr\'{e}chet differentiable at every
point except at $y_{2}$ and $y_{1}$ respectively. Also, there are constants $%
\lambda \in (0,1]$ and $\Lambda \in \lbrack 1,\infty )$ so that 
\begin{equation}
\lambda \Vert y_{1}-y_{2}\Vert _{\mathcal{Y}}\leq M(y_{1},y_{2})\leq \Lambda
\Vert y_{1}-y_{2}\Vert _{\mathcal{Y}},  \label{H2-1}
\end{equation}%
for all $y_{1},y_{2}\in \mathcal{Y}$, and 
\begin{equation}
\lambda \leq ||\nabla M(\cdot ,y_{2})\Vert _{\mathcal{Y}^{\ast }},\quad 
\text{ }\Vert \nabla M(y_{1},\cdot )\Vert _{\mathcal{Y}^{\ast }}\leq \Lambda
,  \label{H2bis}
\end{equation}%
whenever the quantities on the center are defined.

\bigskip

Next, we state and prove our main tool:

\begin{theorem}
\label{Th. locally lipschitz}Let the problem (\ref{problem}), with $U$ open,
and let $w:U\rightarrow \mathbb{R}$ be a continuous, locally bounded
function. In addition to (H1)-(H2), suppose that the following holds:\newline
\ \ (i) the directional derivative $D_{d}f_{x}\left( \cdot \right) $ exists
in $U$ for all unit $d\in \mathcal{Y}$ and all $x\in X.$ Moreover, there is
a constant $C_{0}>0$ such that for all $u\in U$ and all unit $d$ 
\begin{equation}
\inf_{x\in S\left( u\right) }D_{d}f_{x}\left( u\right) >-C_{0};
\label{assumption bdd inf}
\end{equation}%
\newline
(ii) for all unit $d\in \mathcal{Y},$ $w$\ is a viscosity subsolution of\ (%
\ref{eq})\ in $U.$\newline
Then $w$ is locally Lipschitz in $U$.
\end{theorem}

\textit{Proof. }Let $u_{0}\in U$ and take $\eta >0$ so that $w$ is bounded
on $B(u_{0},2\eta )\subset U$. For 
\begin{equation*}
\gamma :=\frac{\lambda }{4\Lambda }<1,
\end{equation*}%
choose any $u_{1},u_{2}\in B(u_{0},\gamma \eta )\subset B(u_{0},\eta )$.
Observe that 
\begin{equation*}
M\left( u_{1},u_{2}\right) \leq \Lambda \Vert u_{1}-u_{2}\Vert \leq 2\Lambda
\gamma \eta =\frac{\lambda }{2}\eta .
\end{equation*}%
Let $C_{1}$ be a positive constant so that 
\begin{equation}
C_{1}\lambda -1-C_{0}>0,  \label{C1}
\end{equation}%
where $C_{0}$ is from \eqref{assumption bdd inf}. Let $h:[0,\infty
)\rightarrow \lbrack 0,\infty )$ be a $C^{1}$ function so that%
\begin{equation}
h^{\prime }(r)\geq C_{1}\text{ \textnormal{for all }}r>0,\quad h(r)=C_{1}r\,%
\text{\textnormal{in}}\,\left[ 0,\frac{\lambda \eta }{2}\right] ,\,\,\text{%
\textnormal{and}}\,h\left( \lambda \eta \right) \geq 2\Vert w\Vert
_{C(B(u_{2},\eta ))}+1.  \label{hh}
\end{equation}%
To get an appropriate test function, define: 
\begin{equation}
g(u):=w(u_{2})+h\left( M(u,u_{2})\right) .  \label{def g}
\end{equation}%
Then $g$ is Fr\'{e}chet differentiable in $U\setminus \left\{ u_{2}\right\} $
and $\left( w-g\right) \left( u_{2}\right) =0$. We will prove that 
\begin{equation}
(w-g)(u)\leq 0\quad \text{for all }u\in B(u_{2},\eta ).  \label{v - phi}
\end{equation}%
Arguing by contradiction, suppose that there is $\tilde{u}\in B\left(
u_{2},\eta \right) $ so that 
\begin{equation}
w(\tilde{u})>g(\tilde{u}).  \label{utilde}
\end{equation}%
Observe that $w-g$ is bounded above and continuous in $B\left( u_{2},\eta
\right) $, then by Theorem \ref{RNP set max}, for every $0<\delta <1$, there
is $h_{\delta }\in \mathcal{Y}^{\ast }$ so that $\Vert h_{\delta }\Vert _{%
\mathcal{Y}^{\ast }}\leq \delta $ and $w-g-h_{\delta }$ attains a maximum in
the closed ball $B\left( u_{2},\eta \right) $. Let $u_{\delta }^{\ast }\in
B\left( u_{2},\eta \right) $ be where $w-g-h_{\delta }$ attains its maximum
over $B\left( u_{2},\eta \right) $ and observe that for $\delta $ small
enough 
\begin{equation}
(w-g-h_{\delta })(u_{\delta }^{\ast })\geq (w-g-h_{\delta })(\tilde{u})>0.
\label{positive}
\end{equation}%
Thus, if $u_{\delta _{i}}^{\ast }=u_{2}$ for a sequence $\delta
_{i}\rightarrow 0^{+}$ as $i\rightarrow \infty $, we obtain from (\ref%
{utilde}) that 
\begin{equation*}
0=(w-g)(u_{2})\geq (w-g)(\tilde{u})>0,
\end{equation*}%
which is a contradiction. Hence, there is some $0<\delta _{0}<1$ so that $%
u_{\delta }^{\ast }\neq u_{2}$ for any $\delta \in (0,\delta _{0})$. From
now on, we assume $\delta <\delta _{0,}$ and prove that $u_{\delta }^{\ast }$
belongs to the interior of the ball $B(u_{2},\eta )$. Indeed, if $u_{\delta
}^{\ast }\in \partial B\left( u_{2},\eta \right) $, then by (\ref{H2-1}), (%
\ref{def g}), and the choice of $h$ \eqref{hh}, we derive 
\begin{equation*}
g(u_{\delta }^{\ast })\geq w(u_{2})+h(\lambda \eta )\geq ||w||_{C\left(
B\left( u_{2},\eta \right) \right) }+1>w(u_{\delta }^{\ast }).
\end{equation*}%
This contradicts (\ref{positive}) for all small enough $\delta $. Thus, $%
u_{\delta }^{\ast }$ is in the interior of the ball $B\left( u_{2},\eta
\right) $. Since $u_{2}\neq u_{\delta }^{\ast }$ and $w$ is a viscosity
subsolution of equation \eqref{eq}, the following holds for any unit $d$ and
for the test function $g+h_{\delta }$%
\begin{equation}
\left\langle -h^{\prime }\left( M(u_{\delta }^{\ast },u_{2})\right) \nabla
M(\cdot ,u_{2})(u_{\delta }^{\ast })-h_{\delta },d\right\rangle +\inf_{x\in
S(u_{\delta }^{\ast })}D_{d}f_{x}(u_{\delta }^{\ast })\leq 0.
\label{eq1xx12}
\end{equation}%
By Hahn-Banach theorem (\cite[Corollary 1.3]{Br}), there is $d_{0}\in 
\mathcal{Y}^{\ast \ast }=\mathcal{Y}$ such that 
\begin{equation}
\left\langle d_{0},\nabla M(\cdot ,u_{2})(u_{\delta }^{\ast })\right\rangle
=\Vert \nabla M(\cdot ,u_{2})(u_{\delta }^{\ast })\Vert _{\mathcal{Y^{\ast }}%
}^{2}  \label{HahnB}
\end{equation}%
and 
\begin{equation*}
\Vert d_{0}\Vert _{\mathcal{Y}}=\Vert \nabla M(\cdot ,u_{2})(u_{\delta
}^{\ast })\Vert _{\mathcal{Y^{\ast }}}.
\end{equation*}%
Plugging the unit vector 
\begin{equation*}
d=-\frac{1}{\Vert \nabla M(\cdot ,u_{2})(u_{\delta }^{\ast })\Vert _{%
\mathcal{Y}^{\ast }}}d_{0}
\end{equation*}%
into \eqref{eq1xx12}, and recalling \eqref{H2bis}, the properties of $h,$
and \eqref{HahnB}, we obtain 
\begin{equation*}
C_{1}\lambda -\left\langle h_{\delta },d\right\rangle +\inf_{x\in
S(u_{\delta }^{\ast })}D_{d}f(x,u_{\delta }^{\ast })\leq 0.
\end{equation*}%
Since $|\left\langle h_{\delta },d\right\rangle |\leq 1$ we have 
\begin{equation*}
C_{1}\lambda -1-C_{0}\leq 0,
\end{equation*}%
which is a contradiction with the choice of $C_{1}$ in \eqref{C1}. In this
way, we have proved \eqref{v - phi}. \newline
To end the proof of this theorem, observe that by the choice of $\gamma $,
it follows that 
\begin{equation*}
\Vert u_{1}-u_{2}\Vert \leq 2\gamma \eta =\frac{\lambda }{2\Lambda }\eta
\leq \frac{\lambda }{2}\eta
\end{equation*}%
and consequently $u_{1}\in B(u_{2},\eta )$ and $h(\Lambda \Vert
u_{1}-u_{2}\Vert )=C_{1}\Lambda \Vert u_{1}-u_{2}\Vert $. Hence by (\ref{def
g}) and \eqref{v - phi}, we get 
\begin{equation*}
w(u_{1})-w(u_{2})\leq h(M(u_{1},u_{2}))\leq h(\Lambda \Vert u_{1}-u_{2}\Vert
)=C_{1}\Lambda \Vert u_{1}-u_{2}\Vert .
\end{equation*}%
Since $u_{1}$ and $u_{2}$ are arbitrary in $B(u_{0},\gamma \eta )$ and $%
u_{0} $ is any point of $U$, we derive the locally Lipschitz regularity of $%
w $ in $U$.\textit{\ \hfill }$\square $\bigskip

Notice that the Lipschitz constant, $C_{1}\Lambda $, that we have found in
the previous proof is the same for all $u_{0}\in U$. Hence, it is easy to
modify this proof to get the following global property:

\begin{theorem}
\label{Th. global Lipschitz} Under the same assumptions as in Theorem \ref%
{Th. locally lipschitz}, if $w$ is bounded and $U=\mathcal{Y}$, then $w$ is
Lipschitz in $\mathcal{Y}.$
\end{theorem}

\begin{remark}
\label{assump H3} We now discuss the feasibility of assumption (H2). This
assumption clearly holds when the underlying norm is Fr\'{e}chet
differentiable. Moreover, observe that by \cite[Theorem 8.19]{FHM}, if $%
\mathcal{Y}^{\ast }$ is separable, then $\mathcal{Y}$ admits an equivalent
norm $\Vert \cdot \Vert _{e}$ which is Fr\'{e}chet differentiable. Hence, we
may define 
\begin{equation*}
M(y_{1},y_{2}):=\Vert y_{1}-y_{2}\Vert _{e},
\end{equation*}%
and observe that \eqref{H2-1} is satisfied for some $\lambda ,\Lambda >0$ so
that%
\begin{equation}
\lambda \Vert y\Vert _{\mathcal{Y}}\leq \Vert y\Vert _{e}\leq \Lambda \Vert
y\Vert _{\mathcal{Y}},\quad y\in \mathcal{Y}.  \label{acotaNorm}
\end{equation}%
Now, we have that $y_{1}\rightarrow M(y_{1},y_{2})$ (resp. $y_{2}\rightarrow
M(y_{1},y_{2})$) is Fr\'{e}chet differentiable at any $y_{1}\neq y_{2}$ ($%
y_{2}\neq y_{1}$) with respect to $\Vert \cdot \Vert _{e}$. Next, we prove
that these mappings are also Fr\'{e}chet differentiable with respect to the
original norm $\Vert \cdot \Vert _{\mathcal{Y}}$ and that in fact their Fr%
\'{e}chet derivatives coincide. Indeed, denoting by $\nabla M(\cdot
,y_{2})(y_{1})$ the Fr\'{e}chet derivative of $y_{1}\rightarrow
M(y_{1},y_{2})$ with respect to $\Vert \cdot \Vert _{e}$, we get for all $%
y\neq 0$, 
\begin{eqnarray*}
0 &\leq &\dfrac{|M(y_{1}+y,y_{2})-M(y_{1},y_{2})-\nabla M(\cdot
,y_{2})(y_{1})y|}{\Vert y\Vert _{\mathcal{Y}}} \\
&\leq &\Lambda \dfrac{|M(y_{1}+y,y_{2})-M(y_{1},y_{2})-\nabla M(\cdot
,y_{2})(y_{1})y|}{\Vert y\Vert _{e}}=o(1),
\end{eqnarray*}%
where we have used \eqref{acotaNorm}. Furthermore, by \cite[pag. 242]{FHM},
we have 
\begin{equation*}
\sup_{y\neq 0}\frac{|\nabla M(\cdot ,y_{2})(y)|}{||y||_{e}}=1.
\end{equation*}%
Another application of \eqref{acotaNorm} gives, for all $y\neq 0,$ that 
\begin{equation*}
\lambda \frac{|\nabla M(\cdot ,y_{2})(y)|}{||y||_{e}}\leq \frac{|\nabla
M(\cdot ,y_{2})(y)|}{||y||}\leq \Lambda \frac{|\nabla M(\cdot ,y_{2})(y)|}{%
||y||_{e}}.
\end{equation*}
The same argument applies to $y_{2}\rightarrow M(y_{1},y_{2})$. Therefore, %
\eqref{H2bis} holds. \bigskip
\end{remark}

\section{The optimal value function as a viscosity solution}

For considering regularity results on $v$, first we state possible
conditions on $f$ and on $\Phi $.

\textbf{Set of assumptions on the function }$f:$

(A1) $\ f$ is continuous in $X\times U$ and satisfies the inf-compactness
property on $gph\Phi .$

(A2) \ For $u\in U$, there is a unit $d(u)$ in $\mathcal{Y}$ so that for all 
$x\in \Phi (u)$, the function $f_{x}(\cdot )=f(x,\cdot )$ is directionally
differentiable at $u$ in the direction $d(u),$ and $\inf_{x\in S\left(
u\right) }D_{d\left( u\right) }f_{x}\left( u\right) >-\infty .$

(A3) \ If $\left\{ x_{k}\right\} \subset \Phi (u)$ converges to $x\in X$,
then 
\begin{equation*}
D_{d(u)}f_{x}(u)\leq \limsup_{k\rightarrow \infty }D_{d(u)}f_{x_{k}}(u),
\end{equation*}%
for $u\in U$, and where $d(u)$ is as in (A2).

(A4) \ For $u\in U$, there is a constant $C_{0}>0$ and a neighborhood $%
\mathcal{N}$ of $u$ such that, for all $u_{1}\in \mathcal{N}$ and all unit $%
d\in \mathcal{Y},$%
\begin{equation}
\inf_{x\in S\left( u_{1}\right) }D_{d}f_{x}\left( u_{1}\right) >-C_{0}.
\label{inf Df}
\end{equation}%
\medskip

Notice that (A1) implies that the optimal set $S\left( u\right) $ is non
empty and compact for any $u\in U$.

For short, we will synthesize the first three conditions on the objective
function $f$ in the following way:\medskip

\textbf{Condition (A):}

(A1) holds in general, and (A2)-(A3) for all unit $d\in \mathbb{R}^{m}$ and
all $u\in U.$ \medskip

Sometimes we will only need that (A2) and (A3) hold true at a particular
point $u_{0}\in U$, and in this case we will write "Condition (A) at $u_{0}$%
". Similarly, "Condition (A) at $u_{0}$ and $d(u_{0})$" means that (A2) and
(A3) are only required to be valid at $u_{0}\in U$ and a unit $d(u_{0})\in 
\mathcal{Y}.$

\medskip

On the side of the feasible set mapping we will analyze two cases. One for a
fixed feasible set 
\begin{equation*}
\emptyset \neq \Phi (u)=\Phi \subset X,
\end{equation*}%
for all $u\in U,$ in a similar fashion as Proposition 4.12 in \cite{BS}. The
other case considers perturbed feasible sets defined by abstract constraints
under the following general assumptions ($X\subset \mathfrak{X}$ and $%
U\subset \mathcal{Y}$ are as usual).\medskip

\textbf{Set of assumptions for a non constant feasible set mapping }$\Phi $:

(B1) $\ \Phi :U\rightarrow 2^{X}$ is defined by 
\begin{equation*}
\Phi (u)=\left\{ x\in X:G(x,u)\in K\right\} ,
\end{equation*}%
where $G:X\times U\rightarrow \mathcal{Z}$ is a continuous function, $%
\mathcal{Z}$ is a given topological vector space, and $K\subset \mathcal{Z}$
is a closed set with non-empty interior.

(B2) \ (Slater-like condition) For $u\in U,$ it holds that $G(x,u)\in \text{%
int }K$ for any $x\in S(u)$. \bigskip

Observe that (B1) implies that $\Phi $ is a closed mapping.

Again, we will synthesize these cases and conditions in the following
way:\medskip

\textbf{Condition (B):}

In the case of a fixed feasible set, it is always closed and, in the case of
perturbed feasible sets, (B1) holds in general and (B2) for all $u\in U$%
.\medskip

As before, we say "Condition (B) at $u_{0}$" when (B2) is only required to
hold at a particular point $u_{0}\in U$.

\begin{remark}[Inf-compactness property and lsc continuity of $v$]
Observe that, under the continuity property of $f$ and, assuming that $\Phi $
is closed, the fact that $f$ satisfies the inf-compactness property at some $%
u_{0}$ guaranties the lsc continuity of $v$ at that point. Indeed, if there
exist a real number $\alpha $ and a compact set $C$ such that 
\begin{equation*}
\emptyset \neq \left\{ x\in \Phi (u):f\left( x,u\right) \leq \alpha \right\}
\subset C
\end{equation*}%
for all $u$\ in some neighborhood of $u_{0}$, then for any sequence $\left\{
u_{k}\right\} $ converging to $u_{0}$ with 
\begin{equation*}
\lim \inf_{u\rightarrow u_{0}}v\left( u\right) =\lim_{k\rightarrow \infty
}v\left( u_{k}\right)
\end{equation*}%
we may find, for $k$ large enough, a sequence $x_{k}\in S\left( u_{k}\right)
\subset C,$ which we assume without loss of generality (w.l.o.g.) that
converges to some $x_{0}\in C.$ So the closedness of $\Phi $ gives that $%
x_{0}\in \Phi (u_{0})$. Hence 
\begin{equation*}
v\left( u_{0}\right) \leq f\left( x_{0},u_{0}\right) =\lim_{k\rightarrow
\infty }f\left( x_{k},u_{k}\right) =\lim_{k\rightarrow \infty }v\left(
u_{k}\right) =\lim \inf_{u\rightarrow u_{0}}v\left( u\right) .
\end{equation*}%
\ \ 
\end{remark}

\begin{remark}[Slater-like condition and usc continuity of $v$]
Suppose that for any $x_{0}\in S(u_{0})$, $G(x_{0},u_{0})\in \text{int }K,$
which is a Slater-like condition. Then it follows immediately the usc
continuity of $v$ at $u_{0}$. In fact, from the continuity of the given
functions $f$ and $G,$ we may find, for any $\varepsilon >0$, neighborhoods $%
V$ and $W$ of $x_{0}\in S(u_{0})$ and $u_{0},$\ respectively, such that 
\begin{equation*}
G\left( x,u\right) \in K\text{ \ \ and \ }f\left( x,u\right) <f\left(
x_{0},u_{0}\right) +\varepsilon =v\left( u_{0}\right) +\varepsilon ,
\end{equation*}%
for all $x\in V$\ and $u\in W.$\ Notice that, for $u\in W,$ $x_{0}\in \Phi
(u)$\ and hence $\Phi \left( u\right) \neq \emptyset $. Thus 
\begin{equation*}
v\left( u\right) \leq v\left( u_{0}\right) +\varepsilon ,
\end{equation*}%
in $W,$ giving the usc continuity of $v$ at $u_{0}.$
\end{remark}

The next lemma states that $v$ is a viscosity solution of a first-order
partial differential equation. This fact will allow us to apply Theorem \ref%
{Th. locally lipschitz} to address the Lipschitz continuity of $v$.

\begin{lemma}
\label{Lema fixed Phi Opt. val.Visc. sol.}For the problem (\ref{problem})
with $U$ open, let $u_{0}\in U$. Assume that $f$ satisfies Conditions (A) at 
$u_{0}$ and $d(u_{0}),$ and (B) at $u_{0}$. Then, the optimal value function 
$v$ is a viscosity solution of the equation%
\begin{equation}
-\left\langle \nabla w\left( u\right) ,d\left( u\right) \right\rangle
+\inf_{x\in S(u)}D_{d(u)}f_{x}(u)=0\,  \label{eq1x}
\end{equation}%
at $u_{0}$.
\end{lemma}

\textit{Proof. }First, consider the case of a fixed closed feasible set. In
view of Theorem \ref{Th. v-continuous} and assumption (A1) of $f,$ we obtain
that $v$ is a continuous function. Let $\eta \in \mathcal{J}^{+}v(u_{0}).$
Take any decreasing sequence, $s_{k}\downarrow 0$ and set $%
u_{k}:=u_{0}+s_{k}d\left( u_{0}\right) .$\ Now, consider optimal solutions $%
x_{k}\in S\left( u_{k}\right) .$ By the inf-compactness property, the fact
that $\Phi $ is closed, and the continuity of $v$, $x_{k}$\ converges to
some $x_{0}\in S\left( u_{0}\right) $ (by passing to a subsequence if
necessary). \newline
Since $\eta \in \mathcal{J}^{+}v(u_{0})$, by considering that $%
f(x_{k},u_{0})\geq f(x_{0},u_{0})=v\left( u_{0}\right) ,$ it follows that 
\begin{eqnarray*}
s_{k}D_{d\left( u_{0}\right) }f_{x_{k}}\left( u_{0}\right) &=&f\left(
x_{k},u_{k}\right) -f\left( x_{k},u_{0}\right) +o\left( s_{k}\right) \\
&\leq &v\left( u_{k}\right) -v\left( u_{0}\right) +o\left( s_{k}\right) \\
&\leq &s_{k}\left\langle \eta ,d\left( u_{0}\right) \right\rangle +o\left(
s_{k}\right) .
\end{eqnarray*}%
Dividing by $s_{k},$ taking $k\rightarrow \infty $, and making use of the
assumption (A3) of $f,$ we obtain 
\begin{equation*}
-\left\langle \eta ,d\left( u_{0}\right) \right\rangle +\inf_{x\in
S(u_{0})}D_{d(u_{0})}f_{x}(u_{0})\leq 0.
\end{equation*}%
Therefore, $v$ is a viscosity subsolution to \eqref{eq1x} at $u_{0}$.\newline
Next, in order to show that $v$ is also a viscosity supersolution to %
\eqref{eq1x} at $u_{0}$, assume that $\eta \in \mathcal{J}^{-}v(u_{0}).$
Take any $x\in S\left( u_{0}\right) ,$ then%
\begin{eqnarray*}
s_{k}D_{d\left( u_{0}\right) }f_{x}\left( u_{0}\right)
&=&f(x,u_{k})-f(x,u_{0})+o\left( s_{k}\right) \\
&\geq &v\left( u_{k}\right) -v\left( u_{0}\right) +o\left( s_{k}\right) \\
&\geq &s_{k}\left\langle \eta ,d\left( u_{0}\right) \right\rangle +o\left(
s_{k}\right) .
\end{eqnarray*}%
Once again, divide by $s_{k}$\ and let $k\rightarrow \infty $, to get $%
D_{d(u_{0})}f_{x}(u_{0})\geq \left\langle \eta ,d\left( u_{0}\right)
\right\rangle $ for all $x\in S\left( u_{0}\right) .$ Hence 
\begin{equation*}
-\left\langle \eta ,d\left( u_{0}\right) \right\rangle +\inf_{x\in
S(u_{0})}D_{d(u_{0})}f_{x}(u_{0})\geq 0,
\end{equation*}%
as desired.\newline
Now, consider the case of perturbed feasible sets satisfying Condition (B)
at $u_{0}.$ In view of the two remarks above, the optimal value function $v$
is continuous at $u_{0}.$ As before, let $\left\{ s_{k}\right\} $ be a
decreasing sequence, $s_{k}\downarrow 0$, and put $%
u_{k}=u_{0}+s_{k}d(u_{0}). $ Taking into account the continuity of $f$ and $%
G,$\ the inf-compactness property, and that $\Phi $ is a closed mapping, we
can obtain, w.l.o.g, optimal solutions $x_{k}\in S\left( u_{k}\right) $
converging to some $x_{0}\in S(u_{0})$. From (B2), there are open
neighborhoods $V$ and $W$ of $x_{0}$ and $u_{0},$ respectively, so that $%
G(x,u)\in \text{int }K,$ i.e. $x\in \Phi \left( u\right) ,$ for all $x\in
V,u\in W.$ In particular, for $k$ large enough, $x_{k}\in \Phi \left(
u_{0}\right) .$ Therefore, we can follow similar steps as above in the case
of fixed feasible sets to get that $v$ is a viscosity subsolution to (\ref%
{eq1x}) at $u_{0}$. \newline
To prove that $v$ is a viscosity supersolution to (\ref{eq1x}) at $u_{0}$,
take any $x$ $\in S\left( u_{0}\right) .$ Once again, (B2) provides open
neighborhoods $V$ and $W$ of $x$ and $u_{0},$ respectively, so that $%
G(x^{\prime },u)\in \text{int }K,$ i.e. $x^{\prime }\in \Phi \left( u\right)
,$ for all $x^{\prime }\in V,u\in W.$ In particular, for $k$ large enough, $%
x $ $\in \Phi \left( u_{k}\right) .$ Once more, we obtain that $v$ is a
viscosity supersolution to (\ref{eq1x}) at $u_{0}$ following the same steps
as in the first case of fixed feasible sets. $\hfill \square $\bigskip

\begin{corollary}
\label{v diffe grad v}Assume all the above conditions and that $\dim 
\mathcal{Y}<\infty .$ If $v$ is differentiable at $u_{0}$ and $f_{x}$ is
differentiable with respect to $u$ at $u_{0}$ for all $x\in S(u_{0}),$ then $%
D_{d}f_{x}(u_{0})$ is constant on $S(u_{0})$ and 
\begin{equation*}
\nabla v(u_{0})=\nabla _{u}f\left( x,u_{0}\right)
\end{equation*}%
for any $x\in S(u_{0})$.
\end{corollary}

\textit{Proof.} From Lemma \ref{Lema fixed Phi Opt. val.Visc. sol.} it
follows that 
\begin{equation*}
-\nabla v\left( u_{0}\right) \cdot d+\inf_{x\in
S(u_{0})}D_{d}f_{x}(u_{0})=0\,
\end{equation*}%
for all unit $d\in \mathcal{Y}$, because $\nabla v(u_{0})\in \mathcal{J}%
^{+}v(u_{0})\cap \mathcal{J}^{-}v(u_{0})$. Furthermore, taking into account
that $D_{d}f_{x}(u_{0})=\nabla _{u}f\left( x,u_{0}\right) \cdot d,$ we have 
\begin{eqnarray*}
\inf_{x\in S(u_{0})}\nabla _{u}f\left( x,u_{0}\right) \cdot d &=&\nabla
v\left( u_{0}\right) \cdot d \\
&=&-\nabla v\left( u_{0}\right) \cdot \left( -d\right) \\
&=&-\inf_{x\in S(u_{0})}\nabla _{u}f\left( x,u_{0}\right) \cdot \left(
-d\right) \\
&=&\sup_{x\in S(u_{0})}\nabla _{u}f\left( x,u_{0}\right) \cdot d.
\end{eqnarray*}%
Hence, $D_{d}f_{x}(u_{0})$ is constant on $S(u_{0})$ and $\nabla
v(u_{0})=\nabla _{u}f\left( x,u_{0}\right) $ for any $x\in S(u_{0})$ as
desired. $\hfill \square $ \bigskip

\begin{theorem}
\label{Th. v loc lipschitz fixed} Let the optimization problem (\ref{problem}%
) where $U$ is open in $\mathcal{Y}.$ Assume hypotheses (H1)-(H2) on $%
\mathcal{Y}$, Conditions (A) and (B), and (A4) for all $u\in U$. Then $v$ is
locally Lipschitz in $U$.
\end{theorem}

\textit{Proof.} It is a straightforward consequence of Theorem \ref{Th.
locally lipschitz} and Lemma \ref{Lema fixed Phi Opt. val.Visc. sol.}. $%
\hfill \square $ \bigskip

\begin{corollary}
Assume all the above conditions. If $\dim \mathcal{Y}<\infty ,$ then $v$ is
differentiable a.e. in $U$ and,$\,$for all unit $d\in \mathcal{Y},$ 
\begin{equation*}
\nabla v(u)\cdot d=\inf_{x\in S(u)}D_{d}f_{x}(u),\text{ \quad a.e. in }U%
\text{.}
\end{equation*}
\end{corollary}

\textit{Proof.} It follows immediately from Lemma \ref{Lema fixed Phi Opt.
val.Visc. sol.} whenever $v$ is differentiable at $u,$ because $\nabla
v(u)\in \mathcal{J}^{+}v(u)\cap \mathcal{J}^{-}v(u)$. $\hfill \square $
\bigskip

\begin{remark}
In the case of $U=\mathcal{Y}$, and a bounded function $f$, if the constant $%
C_{0}$ works globally in (\ref{inf Df}), we also obtain from Theorem \ref%
{Th. global Lipschitz} and Lemma \ref{Lema fixed Phi Opt. val.Visc. sol.}
that $v$ is Lipschitz in $\mathcal{Y}.$
\end{remark}

\section{The case of finite dimensional spaces}

In this section we restrict our analysis to the euclidean space $\mathcal{Y}=%
\mathbb{R}^{m}$. As usual, we identify $\mathcal{Y}^{\ast }$ with $\mathbb{R}%
^{m},$ and use the notations $\left\vert \cdot \right\vert $ and $p\cdot d$
for the euclidean norm and the inner product, respectively. Also, $V\subset
\subset U$ stands for $V\subset \overline{V}\subset U$.

\subsection{Comparison principle and uniqueness of the viscosity solution}

The following lemma gives a comparison principle which will allow us to
characterize $v$ as the unique viscosity solution of a set of equations.

\begin{lemma}[Comparison Principle]
\label{lemma CompPrinc}Let the problem (\ref{problem}) where $U$ is an open
subset of $\mathbb{R}^{m},$ and let $U_{1}$ be any bounded and open subset
such that $U_{1}\subset \subset U$. Assume:\newline
(i) $f:X\times U\rightarrow \mathbb{R}$ has directional derivatives $%
D_{d}f_{x}\left( u\right) $ for all $x\in X,u\in U,$ and all unit $d\in 
\mathbb{R}^{m};$ \newline
(ii) For all $\overline{u}\in \overline{U_{1}},$ all sequences $\left\{
d_{n}\right\} \subset \mathbb{R}^{m},\left\{ u_{n}\right\} \subset U_{1}$\
with 
\begin{equation*}
\left\vert d_{n}\right\vert =1\quad \ \text{and }\quad \lim_{n\rightarrow
\infty }u_{n}=\overline{u},
\end{equation*}%
contain subsequences, which are not relabel, such that 
\begin{equation}
\inf_{x\in S(u_{n})}D_{d_{n}}f_{x}(u_{n})-\inf_{x\in S(\overline{u}%
)}D_{d_{n}}f_{x}(\overline{u})\rightarrow 0\qquad \text{as }n\rightarrow
\infty ,  \label{behaviour direct deriv}
\end{equation}%
and \newline
(iii) $w_{1}\in C\left( \overline{U_{1}}\right) $ and $w_{2}\in C\left( 
\overline{U_{1}}\right) $ are sub- and supersolutions of 
\begin{equation*}
-\nabla w\left( u\right) \cdot d+\inf_{x\in S(u)}D_{d}f_{x}(u)=0\quad \text{%
in }\,U_{1},
\end{equation*}%
respectively, for all unit $d\in \mathbb{R}^{m}$, and that 
\begin{equation}
w_{1}(u)\leq w_{2}(u)\quad \text{for}\,\,u\in \partial U_{1}.
\label{boundary assump}
\end{equation}%
Then 
\begin{equation}
w_{1}\leq w_{2}\quad \text{in }\,U_{1}.  \label{ine sub and super}
\end{equation}
\end{lemma}

\textit{Proof}: Let $u^{\ast }\notin \overline{U_{1}}$. For $\varepsilon
_{n},\delta >0$, with $\varepsilon _{n}\downarrow 0$ as $n\rightarrow \infty 
$, define $\Psi _{n}:\overline{U_{1}}\times \overline{U_{1}}\rightarrow 
\mathbb{R}$ by 
\begin{equation*}
\Psi _{n}(u_{1},u_{2})=w_{1}(u_{1})-w_{2}(u_{2})-\dfrac{|u_{1}-u_{2}|^{2}}{%
2\varepsilon _{n}}-\delta |u_{1}-u^{\ast }|.
\end{equation*}%
There are points $(u_{1,n},u_{2,n})\in \overline{U_{1}}\times \overline{U_{1}%
}$ (depending on $\varepsilon _{n}$ and $\delta $) so that 
\begin{equation*}
\Psi _{n}(u_{1},u_{2})\leq \Psi _{n}(u_{1,n},u_{2,n})
\end{equation*}%
for all $(u_{1},u_{2})\in \overline{U_{1}}\times \overline{U_{1}}$. In
particular 
\begin{equation}
w_{1}(u_{1})-w_{2}(u_{1})=\Psi _{n}(u_{1},u_{1})+\delta |u_{1}-u^{\ast
}|\leq \Psi (u_{1,n},u_{2,n})+\delta |u_{1}-u^{\ast }|,\quad \text{for all }%
u_{1}\in U_{1}.  \label{main ineq 1}
\end{equation}%
W.l.o.g, we may assume that $(u_{1,n},u_{2,n})\rightarrow \left( \overline{u}%
_{1},\overline{u}_{2}\right) \in \overline{U_{1}}\times \overline{U_{1}}$.
Observe that 
\begin{equation*}
\Psi _{n}(u_{1,n},u_{1,n})\leq \Psi _{n}(u_{1,n},u_{2,n}),
\end{equation*}%
which implies that 
\begin{equation*}
w_{2}(u_{2,n})-w_{2}(u_{1,n})\leq -\dfrac{|u_{1,n}-u_{2,n}|^{2}}{%
2\varepsilon _{n}}.
\end{equation*}%
Thus, the boundedness of the function $w_{2}$ and the fact that $\varepsilon
_{n}\downarrow 0$ give 
\begin{equation}
|u_{1,n}-u_{2,n}|\rightarrow 0\quad \text{ as }n\rightarrow \infty ,
\label{conver u eps}
\end{equation}%
and hence $\overline{u}_{1}=\overline{u}_{2}$. \newline
Now, we shall prove that 
\begin{equation}
\liminf_{n\rightarrow \infty }\Psi _{n}(u_{1,n},u_{2,n})\leq 0.
\label{objective}
\end{equation}%
\newline
Indeed, for each $\varepsilon _{n}$, we have two possibilities:\newline
(i) $(u_{1,n},u_{2,n})\in \partial (U_{1}\times U_{1})$\newline
(ii) $(u_{1,n},u_{2,n})\in U_{1}\times U_{1}.$\newline
Suppose that (i) holds for a subsequence $\varepsilon _{n_{k}}\downarrow 0$.
In this case, we may have that $u_{1,n_{k}}\in \partial U_{1}$, and %
\eqref{boundary assump} gives 
\begin{equation*}
\Psi _{n}(u_{1,n_{k}},u_{2,n_{k}})\leq w_{2}(u_{1,n_{k}})-w_{2}(u_{2,n_{k}}).
\end{equation*}%
Also, it may hold that $u_{2,n_{k}}\in \partial U_{1}$ and in this case we
obtain 
\begin{equation*}
\Psi _{n}(u_{1,n_{k}},u_{2,n_{k}})\leq w_{1}(u_{1,n_{k}})-w_{1}(u_{2,n_{k}}).
\end{equation*}%
In either case, \eqref{objective} follows from \eqref{conver u eps} and the
continuity of $w_{1}$ and $w_{2}$.\newline
Next, suppose that $(u_{1,n},u_{2,n})\in U_{1}\times U_{1}$ for all $%
\varepsilon _{n}$ small enough. Let define the test functions 
\begin{equation*}
g_{1,\varepsilon _{n}}:\quad u_{1}\rightarrow w_{2}(u_{2,n})+\dfrac{%
|u_{1}-u_{2,n}|^{2}}{2\varepsilon _{n}}+\delta |u_{1}-u^{\ast }|,
\end{equation*}%
\begin{equation*}
g_{2,\varepsilon _{n}}:\quad u_{2}\rightarrow w_{1}(u_{1,n})-\dfrac{%
|u_{2}-u_{1,n}|^{2}}{2\varepsilon _{n}}-\delta |u_{1,n}-u^{\ast }|.
\end{equation*}%
Then $w_{1}-g_{1,\varepsilon _{n}}=\Psi _{n}(\cdot ,u_{2,n})$ attains a
maximum at $u_{1,n}$, and $w_{2}-g_{2,\varepsilon _{n}}=-\Psi
_{n}(u_{1,n},\cdot )$ attains a minimum at $u_{2,n}$. Moreover, 
\begin{equation*}
\nabla g_{1,\varepsilon _{n}}(u_{1,n})=\frac{u_{1,n}-u_{2,n}}{\varepsilon
_{n}}+\delta \frac{u_{1,n}-u^{\ast }}{|u_{1,n}-u^{\ast }|},
\end{equation*}%
and 
\begin{equation*}
\nabla g_{2,\varepsilon _{n}}(u_{2,n})=\frac{u_{1,n}-u_{2,n}}{\varepsilon
_{n}}.
\end{equation*}%
Since $w_{1}$ is a subsolution and $w_{2}$ is a supersolution, we obtain for
all unit $d$ that 
\begin{equation}
-\left( \frac{u_{1,n}-u_{2,n}}{\varepsilon _{n}}+\delta \frac{%
u_{1,n}-u^{\ast }}{|u_{1,n}-u^{\ast }|}\right) \cdot d+\inf_{x\in
S(u_{1,n})}D_{d}f_{x}(u_{1,n})\leq 0,  \label{eq 1 m}
\end{equation}%
and 
\begin{equation}
-\frac{u_{1,n}-u_{2,n}}{\varepsilon _{n}}\cdot d+\inf_{x\in
S(u_{2,n})}D_{d}f_{x}(u_{2,n})\geq 0.  \label{eq 2 m}
\end{equation}%
From \eqref{eq 1 m} and \eqref{eq 2 m} we derive 
\begin{equation*}
-\delta \left( u_{1,n}-u^{\ast }\right) \cdot d+\inf_{x\in
S(u_{1,n})}D_{d}f_{x}(u_{1,n})-\inf_{x\in S(u_{2,n})}D_{d}f_{x}(u_{2,n})\leq
0.
\end{equation*}%
Choosing 
\begin{equation*}
d_{n}=-\frac{u_{1,n}-u^{\ast }}{|u_{1,n}-u^{\ast }|},
\end{equation*}%
we obtain 
\begin{equation*}
\delta +\inf_{x\in S(u_{1,n})}D_{d_{n}}f_{x}(u_{1,n})-\inf_{x\in
S(u_{2,n})}D_{d_{n}}f_{x}(u_{2,n})\leq 0.
\end{equation*}%
Letting $n\rightarrow \infty ,$ recalling that $\lim_{n\rightarrow \infty
}u_{1,n}=\lim_{n\rightarrow \infty }u_{2,n}$, and appealing to %
\eqref{behaviour direct deriv}, we obtain the contradiction $\delta \leq 0$.
Therefore, just case (i) above may hold and so \eqref{objective} follows.%
\newline
Finally, taking $\liminf $ in \eqref{main ineq 1} gives 
\begin{equation*}
w_{1}(u_{1})-w_{2}(u_{1})\leq \delta |u_{1}-u^{\ast }|
\end{equation*}%
for all $u_{1}\in U_{1}$, and \eqref{ine sub and super} follows by letting $%
\delta \rightarrow 0$. $\hfill \square $ \bigskip

\begin{theorem}
Let the problem (\ref{problem}) where $U\subset \mathbb{R}^{m}$\ is open,
and assume conditions (A) and (B). Let $U_{1}\subset \subset U$ be any
bounded and open subset. If (ii) in Lemma \ref{lemma CompPrinc} holds, then
the optimal value function $v$ is the only continuous function in $\overline{%
U_{1}}$ which satisfies 
\begin{equation*}
w\left( u\right) =\min_{x\in S\left( u\right) }f\left( x,u\right) \quad \ 
\text{on }\partial U_{1}
\end{equation*}%
and 
\begin{equation*}
-\nabla w\left( u\right) \cdot d+\inf_{x\in S(u)}D_{d}f_{x}(u)=0\text{ }%
\quad \text{in }U_{1}
\end{equation*}%
in the viscosity sense for all unit $d\in \mathbb{R}^{m}$.
\end{theorem}

\textit{Proof.} It follows immediately from Lemmas \ref{Lema fixed Phi Opt.
val.Visc. sol.} and \ref{lemma CompPrinc}. $\hfill \square $\bigskip

\subsection{Conditions for optimality: verification technique}

Now, we obtain necessary and sufficient conditions for optimality by
following the line of the verification method for optimal control problems
summarized in \cite[Chapter III]{BCD}.

We consider the problem (\ref{problem}) and assume all the hypotheses
developed in Section 4 for making $v$ a locally Lipschitz viscosity solution
of the equation 
\begin{equation}
-\nabla w\left( u\right) \cdot d+\inf_{x\in S(u)}D_{d}f_{x}(u)=0\text{ }
\label{eq finite}
\end{equation}%
in $U,$ i.e. $U$ is open and conditions (A), (A4) for all $u\in U,$ and (B)
are valid.

Fix $u_{0}\in U.$ We start with the following necessary condition for
optimality.

\begin{proposition}
\label{suf optimality} If $x_{0}\in S\left( u_{0}\right) ,$ then, the
function $w:U\rightarrow \mathbb{R}$ defined by 
\begin{equation*}
w\left( u\right) =f\left( x_{0},u\right)
\end{equation*}%
is a viscosity subsolution of equation (\ref{eq finite}) at $u_{0}$.
\end{proposition}

\textit{Proof.} From condition (B) there is a neighborhood $\mathcal{N}%
\subset U$ of $u_{0}$ such that if $u\in \mathcal{N}$, then $x_{0}\in \Phi
(u).$ Next, consider 
\begin{equation*}
w(u)=f(x_{0},u),\quad u\in \mathcal{N},
\end{equation*}%
and let $g$ be any $C(\mathcal{N})$ function, differentiable at $u_{0}$,
such that 
\begin{equation*}
\max_{\mathcal{N}}(w-g)=(w-g)(u_{0}).
\end{equation*}%
Then, for $u\in \mathcal{N}$, 
\begin{equation*}
v(u)-g(u)\leq w(u)-g(u)\leq w(u_{0})-g(u_{0})=v(u_{0})-g(u_{0}),
\end{equation*}%
which yields that 
\begin{equation*}
\max_{\mathcal{N}}(v-g)=(v-g)(u_{0}).
\end{equation*}%
Since $v$ is a viscosity subsolution, we have 
\begin{equation*}
-\nabla g(u_{0})\cdot d+\inf_{x\in S(u_{0})}D_{d}f_{x}(u_{0})\leq 0,
\end{equation*}%
and the proof is completed. $\hfill \square \bigskip $

A sufficient condition for optimality is the following:

\begin{proposition}
\label{suf opt} Let $x_{0}\in \Phi \left( u_{0}\right) .$ If there exist a
bounded neighborhood $\mathcal{N}\subset \subset U$ of $u_{0}$ and a
viscosity subsolution of (\ref{eq finite}) in $\mathcal{N},$ $w\in C\left( 
\overline{\mathcal{N}}\right) ,$\ that satisfies 
\begin{equation*}
w\left( u\right) \leq \min_{x\in S\left( u\right) }f\left( x,u\right) \text{
\ on }\partial \mathcal{N}\text{ \quad\ and \quad\ }w\left( u\right) \geq
f\left( x_{0},u\right) \text{ for all }u\in \mathcal{N},
\end{equation*}%
then $x_{0}\in S\left( u_{0}\right) .$
\end{proposition}

\textit{Proof.} By the Comparison Principle (Lemma \ref{lemma CompPrinc}), $%
w\left( u\right) \leq v\left( u\right) $ for all $u\in U.$ In particular,
for $u_{0},$ 
\begin{equation*}
v\left( u_{0}\right) \leq f\left( x_{0},u_{0}\right) \leq w\left(
u_{0}\right) \leq v\left( u_{0}\right) .
\end{equation*}%
Hence, $v\left( u_{0}\right) =f\left( x_{0},u_{0}\right) $ and so $x_{0}\in
S\left( u_{0}\right) .$ $\hfill \square \bigskip $

We also provide a sufficient condition for non-optimality.

\begin{proposition}
\label{suf non opt} Let $x_{0}\in \Phi \left( u_{0}\right) .$ If there exist
a bounded neighborhood $\mathcal{N}\subset \subset U$ of $u_{0}$ and a
viscosity supersolution of (\ref{eq finite}) in $\mathcal{N},$ $w\in C\left( 
\overline{\mathcal{N}}\right) ,$\ that satisfies 
\begin{equation*}
w\left( u\right) \geq \min_{x\in S\left( u\right) }f\left( x,u\right) \text{
\ on }\partial \mathcal{N}\text{ \quad\ and \quad\ }w\left( u_{0}\right)
<f\left( x_{0},u_{0}\right) ,
\end{equation*}%
then $x_{0}\notin S\left( u_{0}\right) .$
\end{proposition}

\textit{Proof.} By the comparison principle, $v\leq w$ in $\mathcal{N}$.
Hence, 
\begin{equation*}
v(u_{0})\leq w(u_{0})<f(x_{0},u_{0}),
\end{equation*}%
which implies that $x_{0}\notin S\left( u_{0}\right) .$ $\hfill \square
\bigskip $

\begin{remark}
Observe that the existence of a function $w\in C\left( \overline{\mathcal{N}}%
\right) $ such that $w\geq v$ on $\partial \mathcal{N}$ and 
\begin{equation}
-\nabla w(u)\cdot d+\inf_{x\in \Phi (u)}D_{d}f_{x}(u)\geq 0\quad \text{in }%
\mathcal{N},  \label{aux ineq}
\end{equation}%
in the viscosity sense, implies, by Proposition \ref{suf non opt}, that $%
x_{0}\notin S(u_{0})$. The inequality \eqref{aux ineq} does not require to
find $S(u)$.
\end{remark}

\subsection{Generalized derivative and differentiability properties of $v$}

Here we utilize the viscosity solution approach to consider generalized
solutions, in Clarke sense, of the optimal value function $v$, for $\mathcal{%
Y}=\mathbb{R}^{m}$. All the definitions and properties stated with no proof
can be found in \cite{Clk}.

\begin{defn}
Let $u_{0}\in U$ and let $w:U\rightarrow \mathbb{R}$ be Lipschitz in a
neighborhood of $u_{0}$. The generalized derivative $D^{0}w(u_{0}):\mathbb{R}%
^{m}\rightarrow \mathbb{R}$ is defined by 
\begin{equation*}
D^{0}w(u_{0})y:=\limsup_{u\rightarrow u_{0},\,h\downarrow 0}\dfrac{%
w(u+hy)-w(u)}{h},
\end{equation*}%
for $y\in \mathbb{R}^{m}$.\newline
The generalized gradient $\partial w(u_{0})$ of $w$ at $u_{0}$ is given by 
\begin{equation*}
\partial w(u_{0})=\left\{ p\in \mathbb{R}^{m}:p\cdot y\leq D^{0}w(u_{0})y%
\text{ for all }y\in \mathbb{R}^{m}\right\} .
\end{equation*}
\end{defn}

The generalized gradient coincides with the derivative when the function $w$
is strictly differentiable, e.g. when $w$ is $C^{1}$. If $w$ is convex on $U$%
, then $\partial w(u)$ is the subdifferential at $u$ in the sense of convex
analysis, and $D^{0}w(u)y$ is the usual directional derivative $D_{y}w(u),$\
see \cite[Proposition 2.2.7]{Clk}.

By \cite[Proposition 2.1.5]{Clk} the set valued mapping $\partial
w:U\rightrightarrows \mathbb{R}^{m}$ is a closed and usc mapping. Moreover,
we have the next result related to super- and subgradients (see \cite[%
Theorem 1.4]{F}).

\begin{proposition}
\label{rela subdifferential generalized}If $w$ is Lipschitz in a
neighborhood of $u_{0}$, then 
\begin{equation*}
\mathcal{J}^{+}w(u_{0})\cup \mathcal{J}^{-}w(u_{0})\subset \partial w(u_{0}).
\end{equation*}
\end{proposition}

An alternative definition of the generalized gradient is given by the
following proposition (see \cite[Theorem 2.5.1]{Clk}).

\begin{proposition}
\label{characterizetion generalized} If $w$ is Lipschitz in a neighborhood
of $u_{0}$, then 
\begin{equation*}
\partial w(u_{0})=co\left\{ \lim_{l\rightarrow \infty }\nabla
w(u_{l}):u_{l}\rightarrow u_{0}\right\} ,
\end{equation*}%
where the limit is taken over all sequences $\left\{ u_{l}\right\} $
converging to $u_{0}$ so that $\nabla w(u_{l})$ does exist and $\left\{
\nabla w(u_{l})\right\} $ is a converging subsequence.
\end{proposition}

Next, given an operator \ $F:U\times \mathbb{R}\times \mathbb{R}%
^{m}\rightarrow \mathbb{R},$ we introduce the definition of generalized
solutions.

\begin{defn}
A locally Lipschitz function $w:U\rightarrow \mathbb{R}$ is a generalized
solution of 
\begin{equation*}
F(\cdot ,w,\nabla w)=0
\end{equation*}%
at $u_{0}\in U$, if 
\begin{equation*}
\max_{p\in \partial w(u_{0})}F(u_{0},w(u_{0}),p)=0.
\end{equation*}
\end{defn}

We now consider the optimal value function $v$ and provide conditions for $v$
to be a generalized solution of the equation (\ref{eq finite}). For that
recall that Lemma \ref{Lema fixed Phi Opt. val.Visc. sol.} and Theorem \ref%
{Th. v loc lipschitz fixed} give conditions under which $v$ is a locally
Lipschitz viscosity solution of this equation.

\begin{theorem}
\label{viscosity implies generalized} Let the problem (\ref{problem}), where 
$f$ satisfies (A1). Let $u_{0}\in U$ and $d\in \mathbb{R}^{m}.$ Suppose that 
$D_{d}f$ is continuous in $X\times \mathcal{N}$ for a neighborhood $\mathcal{%
N}$ of $u_{0}$, and assume that $w$ is a locally Lipschitz viscosity
solution of 
\begin{equation}
-\nabla w(u)\cdot d+\inf_{x\in S(u)}D_{d}f_{x}(u)=0  \label{equation for v}
\end{equation}%
in $\mathcal{N}$. Moreover, suppose that $\mathcal{J}^{-}w(u_{0})\neq
\emptyset $. Then, $w$ is also a generalized solution of 
\eqref{equation for
v} at $u_{0}$.
\end{theorem}

\textit{Proof.} Since $w$ is a viscosity supersolution, the fact that $%
\mathcal{J}^{-}w(u_{0})\neq \emptyset $ and Proposition \ref{rela
subdifferential generalized}, yield 
\begin{equation}
\max_{p\in \partial w(u_{0})}\left( -p\cdot d+\inf_{x\in
S(u_{0})}D_{d}f_{x}(u_{0})\right) \geq 0.  \label{max gr eq 0}
\end{equation}%
To prove the opposite inequality, take any $p_{0}\in \partial w(u_{0})$. By
Proposition \ref{characterizetion generalized}, we may write 
\begin{equation*}
p_{0}=\sum_{i=1}^{k}\lambda _{i}p_{i},
\end{equation*}%
where $\lambda _{i}\geq 0$ for all $i$, $\sum_{i=1}^{k}\lambda _{i}=1$, and 
\begin{equation*}
p_{i}=\lim_{l\rightarrow \infty }\nabla w(u_{l}^{i}),
\end{equation*}%
with $u_{l}^{i}\rightarrow u_{0}$ as $l\rightarrow \infty ,$ for all $i$.
W.l.o.g., assume that $u_{l}^{i}\in \mathcal{N}$ for all $i$ and $l$. Since $%
w$ is a viscosity solution in $\mathcal{N}$ of \ref{equation for v}, we have 
\begin{equation*}
-\nabla w(u_{l}^{i})\cdot d+\inf_{x\in S(u_{l}^{i})}D_{d}f_{x}(u_{l}^{i})=0.
\end{equation*}%
because $\nabla w(u_{l}^{i})\in \mathcal{J}^{+}w(u_{l}^{i})\cup \mathcal{J}%
^{-}w(u_{l}^{i})$. Now, $S(u_{l}^{i})$ is compact and $D_{d}f_{(\cdot
)}(u_{l}^{i})$ is continuous, hence there is some $x_{l}^{i}\in S(u_{l}^{i})$
so that 
\begin{equation*}
\inf_{x\in S(u_{l}^{i})}D_{d}f_{x}(u_{l}^{i})=D_{d}f_{x_{l}^{i}}(u_{l}^{i}),
\end{equation*}%
and then 
\begin{equation}
-\nabla w(u_{l}^{i})\cdot d+D_{d}f_{x_{l}^{i}}(u_{l}^{i})=0.
\end{equation}
By the inf-compactness property, there exists $x^{i}$ so that $%
x_{l}^{i}\rightarrow x^{i}$, passing to a subsequence if necessary. Observe
that $x^{i}\in S(u_{0})$. Finally, by letting $l\rightarrow \infty $, 
\begin{equation*}
-p_{i}\cdot d+D_{d}f_{x^{i}}(u_{0})=0.
\end{equation*}%
Hence, 
\begin{equation*}
-p_{0}\cdot d+\inf_{x\in S(u_{0})}D_{d}f_{x}(u_{0})\leq
\sum_{i=1}^{k}\lambda _{i}\left( -p_{i}\cdot d+D_{d}f_{x^{i}}(u_{0})\right)
=0.
\end{equation*}%
Thus, 
\begin{equation}
\max_{p\in \partial w(u_{0})}\left( -p\cdot d+\inf_{x\in
S(u_{0})}D_{d}f_{x}(u_{0})\right) \leq 0.  \label{max less eq 0}
\end{equation}%
(\ref{max gr eq 0}) and (\ref{max less eq 0}) complete the proof. \textit{\
\hfill }$\square $

\begin{corollary}
\label{Cor 22} Let the problem (\ref{problem}) that satisfies (A1) and
condition (B). Let $u_{0}\in U$ and suppose the existence of a neighborhood $%
\mathcal{N}$ of it such that $D_{d}f$ is continuous in $X\times \mathcal{N}$
for all unit $d\in \mathbb{R}^{m}$. If $\mathcal{J}^{-}v(u_{0})\neq
\emptyset ,$ then \newline
(i) 
\begin{equation*}
-p\cdot d+\inf_{x\in S(u_{0})}D_{d}f_{x}(u_{0})=0,
\end{equation*}%
for all $p\in \mathcal{J}^{-}v(u_{0})$ and all unit $d\in \mathbb{R}^{m};$ 
\newline
(ii) $\mathcal{J}^{-}v(u_{0})$ is a singleton;\newline
(iii) $D_{d}f(\cdot ,u_{0})$ is constant in $S(u_{0})$ for all unit $d\in 
\mathbb{R}^{m};$ and,\newline
(iv) if $v$ is differentiable at $u_{0}$, then 
\begin{equation}
\nabla v(u_{0})=\nabla _{u}f(x,u_{0}),  \label{main eq derivative}
\end{equation}%
where $x$ is any point in $S(u_{0})$.
\end{corollary}

\textit{Proof}: Notice that the hypotheses imply all the conditions in Lemma %
\ref{Lema fixed Phi Opt. val.Visc. sol.} and Theorem \ref{Th. v loc
lipschitz fixed} that make $v$ a locally Lipschitz viscosity solution of (%
\ref{equation for v}) in $\mathcal{N}.$ Hence, an application of (\ref{max
gr eq 0}) and (\ref{max less eq 0}) in the proof of Theorem \ref{viscosity
implies generalized} gives (i). \newline
(ii) Let $p,q\in \mathcal{J}^{-}v(u_{0}),$ then 
\begin{equation*}
p\cdot d=\inf_{x\in S(u_{0})}D_{d}f_{x}(u_{0})=q\cdot d,
\end{equation*}%
for all unit $d\in \mathbb{R}^{m}$, which yields $p=q.$ \newline
(iii) Here $D_{d}f_{x}(u_{0})=\nabla _{u}f(x,u_{0})\cdot d.$ Now, in view of
(i), we have for $p\in \mathcal{J}^{-}v(u_{0})$ and for all unit $d\in 
\mathbb{R}^{m}$ that 
\begin{eqnarray*}
\inf_{x\in S(u_{0})}\nabla _{u}f(x,u_{0})\cdot d &=&-p\cdot (-d) \\
&=&-\inf_{x\in S(u_{0})}\left\{ \nabla _{u}f(x,u_{0})\cdot (-d)\right\}  \\
&=&\sup_{x\in S(u_{0})}\nabla _{u}f(x,u_{0})\cdot d)
\end{eqnarray*}%
Thus, $D_{d}f(\cdot ,u_{0})$ is constant in $S(u_{0})$.\newline
(iv) Observe that when $v$ is differentiable at $u_{0}$, we have 
\begin{equation*}
\nabla v(u_{0})\cdot d=\nabla _{u}f(x,u_{0})\cdot d
\end{equation*}%
for all unit $d$ and any $x\in S(u_{0})$. Hence, \eqref{main eq
derivative} follows.\textit{\ \hfill }$\square $\bigskip 

Observe that we have already presented in Corollary \ref{v diffe grad v} a
similar result to (iii) and (iv) under slight different conditions.
Moreover, \cite[Remark 4.14]{BS} discusses conditions so that $D_{d}f(\cdot
,u_{0})$ is constant in $S(u_{0})$ when $U$ is any Banach space.

Finally, we include some results for the convex setting. Notice that $v$ is
convex in $U$ whenever the optimization problem (\ref{problem}) is convex,
i.e. $X$ and $U$ are convex sets, $f=f(x,u)$ is convex in $(x,u)\in X\times
U,$ and $gph\Phi \subset X\times U$ is also convex.

\begin{proposition}
\label{prop. convex} Let the problem (\ref{problem}) with $U$ open and let $%
u_{0}\in U.$\ Assume (A1) and condition (B). Suppose the existence of a
neighborhood $\mathcal{N}$ of $u_{0}$ such that $D_{d}f$ is continuous in $%
X\times \mathcal{N}$ for all unit $d\in \mathbb{R}^{m}$. If $v$ is convex
near $u_{0}$, then $v$ is differentiable at $u_{0}$.
\end{proposition}

\textit{Proof}: The convexity of $v$ gives that it is locally Lipschitz and,
moreover, that $\mathcal{J}^{-}v(u_{0})\neq \emptyset .$ By virtue of
Corollary \ref{Cor 22}, $\mathcal{J}^{-}v(u_{0})$ is a singleton, which
yields the differentiability of $v$ at $u_{0}$ (see, e.g., \cite[Theorem 25.1%
]{R}).\textit{\ \hfill }$\square $

\begin{remark}
Proposition 2.2.7 in \cite{Clk} yields that for $v$ convex on $U,$ the
generalized gradient $\partial v(u_{0})$ is the usual subdifferential in
convex analysis and, moreover, $D^{0}v(u_{0})y$\ is the usual directional
derivative $D_{d}v(u_{0}).$\ Proposition \ref{prop. convex} goes further and
provides conditions so that a convex optimal value function is actually
differentiable at the considered point $u_{0}$.
\end{remark}

\bigskip

\begin{theorem}
Assume that both $\mathcal{Y}$ and $\mathfrak{X}$ are finite dimensional
spaces. Let the problem (\ref{problem}), with $X$ and $U$ open sets, and $%
f\in C^{1}(X\times U).$ Suppose that (A1) and condition (B) hold, and let $%
u_{0}\in U$. If $v$ is convex near $u_{0}$, then it is $C^{1}$ in a
neighborhood of $u_{0}$.
\end{theorem}

\textit{Proof. }Let $\mathcal{O}\subset \subset X$ be an open and bounded
neighborhood of $S(u_{0})$, and let $\mathcal{N}\subset \subset U$\textit{\ }%
be an open and bounded neighborhood of $u_{0}$ where $v$ is convex. Fix a
unit $d\in \mathcal{Y}$. Then $D_{d}f(x,u)$ is uniformly continuous in $%
\overline{\mathcal{O}}\times \overline{\mathcal{N}}$. Hence, for fixed $\eta
>0$, there is $\delta =\delta (\eta )>0$ so that 
\begin{equation}
|D_{d}f(x_{1},u_{1})-D_{d}f(x_{2},u_{2})|<\eta 
\label{uniform continuiti D_d}
\end{equation}%
for any $(x_{1},u_{1}),(x_{2},u_{2})\in \overline{\mathcal{O}}\times 
\overline{\mathcal{N}}$ with $\max \left\{
|x_{1}-x_{2}|,|u_{1}-u_{2}|\right\} <\delta $. Moreover, we assume that the $%
\delta -$neighborhood of the compact set $S(u_{0})$ given by 
\begin{equation*}
\mathcal{O}_{\delta }:=\left\{ x\in \mathbb{R}^{n}:\text{ dist }%
(x,S(u_{0}))<\delta \right\} 
\end{equation*}%
is contained in $\mathcal{O}$. Now, $S$ is usc because it is a closed and
compact valued mapping, thus there is a neighborhood $\mathcal{N}_{1}\subset 
\mathcal{N}\cap B\left( u_{0},\delta \right) $ of $u_{0}$ such that 
\begin{equation*}
u\in \mathcal{N}_{1}\quad \text{implies }\,S(u)\subset \mathcal{O}_{\delta }.
\end{equation*}%
In particular, for all $x_{u}\in S(u)$, with $u\in \mathcal{N}_{1}$, there
is $x(u)\in S(u_{0})$ such that 
\begin{equation}
|x_{u}-x(u)|=\text{dist }(x_{u},S(u_{0}))<\delta .  \label{cont trayectory}
\end{equation}%
Also, observe that, by Proposition \ref{prop. convex}, $v$ is differentiable
in a neighborhood $\mathcal{N}_{2}\subset \mathcal{N}_{1}$ of $u_{0}$. Thus,
for all $u\in \mathcal{N}_{2}$ and all unit $d$, it holds that 
\begin{equation*}
\left\vert \nabla v(u)\cdot d-\nabla v(u_{0})\cdot d\right\vert =\left\vert
D_{d}f(x_{u},u)-D_{d}f(x(u),u_{0})\right\vert <\eta ,
\end{equation*}%
where we have used Corollary \ref{Cor 22} (iv), \eqref{cont trayectory}, and %
\eqref{uniform continuiti D_d}. Finally, if $u_{1},u_{2}\in \mathcal{N}_{2},$
then 
\begin{equation*}
\left\vert \nabla v(u_{1})\cdot d-\nabla v(u_{2})\cdot d\right\vert \leq
\left\vert \nabla v(u_{1})\cdot d-\nabla v(u_{0})\cdot d\right\vert
+\left\vert \nabla v(u_{2})\cdot d-\nabla v(u_{0})\cdot d\right\vert <2\eta ,
\end{equation*}%
which completes the proof. \textit{\ \hfill }$\square $\bigskip 


\end{document}